\documentclass{amsart}
\usepackage{latexsym}
\usepackage{amssymb}
\usepackage{amsmath}

\begin{document}


\newtheorem{theorem}{Theorem}
\newtheorem{problem}{Problem}
\newtheorem{definition}{Definition}
\newtheorem{lemma}{Lemma}
\newtheorem{proposition}{Proposition}
\newtheorem{corollary}{Corollary}
\newtheorem{example}{Example}
\newtheorem{conjecture}{Conjecture}
\newtheorem{algorithm}{Algorithm}
\newtheorem{exercise}{Exercise}
\newtheorem{xample}{Example}
\newtheorem{remarkk}{Remark}

\newcommand{\be}{\begin{equation}}
\newcommand{\ee}{\end{equation}}
\newcommand{\bea}{\begin{eqnarray}}
\newcommand{\eea}{\end{eqnarray}}
\newcommand{\beq}[1]{\begin{equation}\label{#1}}
\newcommand{\eeq}{\end{equation}}
\newcommand{\beqn}[1]{\begin{eqnarray}\label{#1}}
\newcommand{\eeqn}{\end{eqnarray}}
\newcommand{\beaa}{\begin{eqnarray*}}
\newcommand{\eeaa}{\end{eqnarray*}}
\newcommand{\req}[1]{(\ref{#1})}

\newcommand{\lip}{\langle}
\newcommand{\rip}{\rangle}
\newcommand{\uu}{\underline}
\newcommand{\oo}{\overline}
\newcommand{\La}{\Lambda}
\newcommand{\la}{\lambda}
\newcommand{\eps}{\varepsilon}
\newcommand{\om}{\omega}
\newcommand{\Om}{\Omega}
\newcommand{\ga}{\gamma}
\newcommand{\ka}{\kappa}
\newcommand{\rrr}{{\Bigr)}}
\newcommand{\qqq}{{\Bigl\|}}

\newcommand{\dint}{\displaystyle\int}
\newcommand{\dsum}{\displaystyle\sum}
\newcommand{\dfr}{\displaystyle\frac}
\newcommand{\bige}{\mbox{\Large\it e}}
\newcommand{\integers}{{\Bbb Z}}
\newcommand{\rationals}{{\Bbb Q}}
\newcommand{\reals}{{\rm I\!R}}
\newcommand{\realsd}{\reals^d}
\newcommand{\realsn}{\reals^n}
\newcommand{\NN}{{\rm I\!N}}
\newcommand{\DD}{{\rm I\!D}}
\newcommand{\degree}{{\scriptscriptstyle \circ }}
\newcommand{\dfn}{\stackrel{\triangle}{=}}
\def\complex{\mathop{\raise .45ex\hbox{${\bf\scriptstyle{|}}$}
     \kern -0.40em {\rm \textstyle{C}}}\nolimits}
\def\hilbert{\mathop{\raise .21ex\hbox{$\bigcirc$}}\kern -1.005em {\rm\textstyle{H}}} 
\newcommand{\RAISE}{{\:\raisebox{.6ex}{$\scriptstyle{>}$}\raisebox{-.3ex}
           {$\scriptstyle{\!\!\!\!\!<}\:$}}} 

\newcommand{\hh}{{\:\raisebox{1.8ex}{$\scriptstyle{\degree}$}\raisebox{.0ex}
           {$\textstyle{\!\!\!\! H}$}}}

\newcommand{\OO}{\won}
\newcommand{\calA}{{\mathcal A}}
\newcommand{\calB}{{\mathcal B}}
\newcommand{\calC}{{\cal C}}
\newcommand{\calD}{{\cal D}}
\newcommand{\calE}{{\cal E}}
\newcommand{\calF}{{\mathcal F}}
\newcommand{\calG}{{\cal G}}
\newcommand{\calH}{{\cal H}}
\newcommand{\calK}{{\cal K}}
\newcommand{\calL}{{\mathcal L}}
\newcommand{\calM}{{\mathcal M}}
\newcommand{\calO}{{\cal O}}
\newcommand{\calP}{{\cal P}}
\newcommand{\calT}{{\mathcal T}} 
\newcommand{\calU}{{\mathcal U}}
\newcommand{\calX}{{\cal X}}
\newcommand{\calY}{{\mathcal Y}}
\newcommand{\calZ}{{\mathcal Z}}
\newcommand{\calXX}{{\cal X\mbox{\raisebox{.3ex}{$\!\!\!\!\!-$}}}}
\newcommand{\calXXX}{{\cal X\!\!\!\!\!-}}
\newcommand{\gi}{{\raisebox{.0ex}{$\scriptscriptstyle{\cal X}$}
\raisebox{.1ex} {$\scriptstyle{\!\!\!\!-}\:$}}}
\newcommand{\intsim}{\int_0^1\!\!\!\!\!\!\!\!\!\sim}
\newcommand{\intsimt}{\int_0^t\!\!\!\!\!\!\!\!\!\sim}
\newcommand{\pp}{{\partial}}
\newcommand{\al}{{\alpha}}
\newcommand{\sB}{{\cal B}}
\newcommand{\sL}{{\cal L}}
\newcommand{\sF}{{\cal F}}
\newcommand{\sE}{{\cal E}}
\newcommand{\sX}{{\cal X}}
\newcommand{\R}{{\rm I\!R}}
\renewcommand{\L}{{\rm I\!L}}
\newcommand{\vp}{\varphi}
\newcommand{\N}{{\rm I\!N}}
\def\ooo{\lip}
\def\ccc{\rip}
\newcommand{\ot}{\hat\otimes}
\newcommand{\rP}{{\Bbb P}}
\newcommand{\bfcdot}{{\mbox{\boldmath$\cdot$}}}

\renewcommand{\varrho}{{\ell}}
\newcommand{\dett}{{\textstyle{\det_2}}}
\newcommand{\sign}{{\mbox{\rm sign}}}
\newcommand{\TE}{{\rm TE}}
\newcommand{\TA}{{\rm TA}}
\newcommand{\E}{{\rm E\,}}
\newcommand{\won}{{\mbox{\bf 1}}}
\newcommand{\Lebn}{{\rm Leb}_n}
\newcommand{\Prob}{{\rm Prob\,}}
\newcommand{\sinc}{{\rm sinc\,}}
\newcommand{\ctg}{{\rm ctg\,}}
\newcommand{\loc}{{\rm loc}}
\newcommand{\trace}{{\,\,\rm trace\,\,}}
\newcommand{\Dom}{{\rm Dom}}
\newcommand{\ifff}{\mbox{\ if and only if\ }}
\newcommand{\nproof}{\noindent {\bf Proof:\ }}
\newcommand{\remark}{\noindent {\bf Remark:\ }}
\newcommand{\remarks}{\noindent {\bf Remarks:\ }}
\newcommand{\note}{\noindent {\bf Note:\ }}

\newcommand{\boldx}{{\bf x}}
\newcommand{\boldX}{{\bf X}}
\newcommand{\boldy}{{\bf y}}
\newcommand{\boldR}{{\bf R}}
\newcommand{\uux}{\uu{x}}
\newcommand{\uuY}{\uu{Y}}

\newcommand{\limn}{\lim_{n \rightarrow \infty}}
\newcommand{\limN}{\lim_{N \rightarrow \infty}}
\newcommand{\limr}{\lim_{r \rightarrow \infty}}
\newcommand{\limd}{\lim_{\delta \rightarrow \infty}}
\newcommand{\limM}{\lim_{M \rightarrow \infty}}
\newcommand{\limsupn}{\limsup_{n \rightarrow \infty}}

\newcommand{\ra}{ \rightarrow }

\newcommand{\ARROW}[1]
  {\begin{array}[t]{c}  \longrightarrow \\[-0.2cm] \textstyle{#1} \end{array} }

\newcommand{\AR}
 {\begin{array}[t]{c}
  \longrightarrow \\[-0.3cm]
  \scriptstyle {n\rightarrow \infty}
  \end{array}}

\newcommand{\pile}[2]
  {\left( \begin{array}{c}  {#1}\\[-0.2cm] {#2} \end{array} \right) }

\newcommand{\floor}[1]{\left\lfloor #1 \right\rfloor}

\newcommand{\mmbox}[1]{\mbox{\scriptsize{#1}}}

\newcommand{\ffrac}[2]
  {\left( \frac{#1}{#2} \right)}

\newcommand{\one}{\frac{1}{n}\:}
\newcommand{\half}{\frac{1}{2}\:}

\def\le{\leq}
\def\ge{\geq}
\def\lt{<}
\def\gt{>}

\def\squarebox#1{\hbox to #1{\hfill\vbox to #1{\vfill}}}
\newcommand{\nqed}{\hspace*{\fill}
          \vbox{\hrule\hbox{\vrule\squarebox{.667em}\vrule}\hrule}\bigskip}

\title{Entropy, Invertibility and Variational Calculus  of the Adapted Shifts on
  Wiener space}

\author{Ali  S\"uleyman  \"Ust\"unel}
\maketitle
\begin{abstract}
\noindent
 In this work we study the necessary and sufficient
conditions for a positive random variable whose expectation under
the Wiener measure is one, to be represented as the Radon-Nikodym
derivative of the image of the Wiener measure under  an adapted
perturbation of identity
 with the help of the associated innovation process. We prove that the innovation conjecture
holds if and only if the original process is almost surely
invertible. We also  give 
variational characterizations of the invertibility of the
perturbations of identity and the representability of a  positive
random variable  whose  total mass is  equal to unity.  We prove in
particular  that an
adapted perturbation of identity $U=I_W+u$ satisfying the Girsanov
theorem, is invertible if and only if the kinetic energy of $u$ is
equal to the entropy of the measure induced with the action of $U$
on the Wiener measure $\mu$, in  other words $U$  is invertible iff
$$
\half \int_W|u|_H^2d\mu=\int_W \frac{dU\mu}{d\mu}\log\frac{dU\mu}{d\mu}d\mu\,.
$$
 The relations with the Monge-Kantorovitch
measure transportation are also studied.  An  application of these
 results  to a variational problem related to large deviations is also
 given.
\end{abstract}

\tableofcontents
\section{Introduction}
This paper is devoted to the study of the following question:
assume that $(W,H,\mu)$ is the classical Wiener space, i.e.,
$W=C_0([0,1],\R^d)$, $H$ is the corresponding Cameron-Martin space
consisting of the absolutely continuous, $\R^d$-valued functions
on $[0,1]$ with square integrable derivatives. Assume that $L$ is
a strictly positive random variable whose expectation with respect
to $\mu$ is one. We suppose that there exits a map $U:W\to W$ of
the form $U=I_W+u$, with $u:W\to H$ such that $\dot{u}$ is adapted
to the filtration of the Wiener space and  that $L$ is represented
by $U$, i.e.
$$
\frac{dU\mu}{d\mu}=L\,.
$$
We suppose  also that
$$
E[\rho(-\delta u)]=1\,,
$$
where
$$
\rho(-\delta
u)=\exp\left[-\int_0^1(\dot{u}_s,dW_s)-\frac{1}{2}\int_0^1|\dot{u}_s|^2ds\right]\,.
$$
Then $U\mu$ is equivalent to $\mu$ and the corresponding
Radon-Nikodym derivative $L$ can be represented as an exponential
martingale $\rho(-\delta v)$ where $v:W\to H$ satisfies similar
properties as those satisfied by $u$.  The question we adress is:
what are the relations satisfied by the couple $(u,v)$? For
instance, if $U$ and $V=I_W+v$ are inverse to each other then the
situation described above happens.  However,  due to the
celebrated example of Tsirelson (cf. \cite{Tsi}), we know that this is not the only
case. We  concentrate ourselves particularly to  this case with
the help of associated innovation processes, in terms of which we
give necessary and sufficient conditions for the representability
(c.f. \cite{FUZ}) of a strictly positive density and for the
invertibility of the associated perturbation of identity. The
innovation approach leads to a nice result which characterizes the
invertibility of an adapted shift in terms of the relative entropy
of the measure which it induces. Namely, assume that $U=I_W+u$ as
above, then it is invertible if and only if the relative entropy
$H(U\mu|\mu)$ is equal to the kinetic energy of $u$, i.e.,
$$
H(U\mu|\mu)=\half E\int_0^1|\dot{u}_s|^2ds\,.
$$
In Physics the notion of entropy is an indication for the
number of accessible states; here it is a  remarkable fact that the
relative entropy behaves as the physical entropy in the sense that if
the system has just  enough kinetic energy to fulfill the accessible states, i.e., if
this energy is equal to the relative entropy of the probability
distribution that it creates then the mapping is invertible. Besides,
in general it is always larger or equal to  the latter. 

We apply this considerations to the innovation problem of the
filtering. Namely it is a celebrated question whether the sigma
algebra generated by the observation process is equal to that of the
innovation process. The case the signal is independent of the noise
has been solved in \cite{A-M}, here we solve this problem in terms of
the entropy of the observed  system.

If we represent a density of the form $L=\rho(-\delta v)$ by
$U=I_W+u$, then, modulo some integrability hypothesis, the
Girsanov theorem implies that $(I_W+v)\circ U=V\circ U$ is a
Wiener process. We study then the properties of $U\circ V$ using
similar techniques. The relations with the Monge transportation
are also exhibited.

In the final part  we use the variational methods to
characterize the invertibility and  representability of densities.
As an  application we give some new results for a
particular case studied in \cite{B-D}. Namely we give an explicit
characterization of the solution of the minimization problem
$$
\inf\left(E\left[f\circ U+\half |u|_H^2\right]\right)\,,
$$
with the help of the entropic characterization of the invertibility
explained above, 
where the inf is taken in the space of adapted, $H$-valued Wiener
functionals with finite energy and $f$ is a $1$-convex Wiener
functional in the Sobolev space $\DD_{2,1}(H)$.
\section{Preliminaries and notation}
\label{preliminaries}
Let $W$ be the classical Wiener  space with  the Wiener
measure $\mu$. The
corresponding Cameron-Martin space is denoted by $H$. Recall that the
injection $H\hookrightarrow W$ is compact and its adjoint is the
natural injection $W^\star\hookrightarrow H^\star\subset
L^2(\mu)$. A subspace $F$ of $H$ is called regular if the
corresponding orthogonal projection
has a continuous extension to $W$, denoted again  by the same letter.
It is well-known that there exists an increasing sequence of regular
subspaces $(F_n,n\geq 1)$, called total,  such that $\cup_nF_n$ is
dense in $H$ and in $W$. Let $\sigma(\pi_{F_n})${\footnote{For the notational
  simplicity, in the sequel we shall denote  it by  $\pi_{n}$.}}  be the
$\sigma$-algebra generated by $\pi_{F_n}$, then  for any  $f\in
L^p(\mu)$, the martingale  sequence
$(E[f|\sigma(\pi_{F_n})],n\geq 1)$
converges to $f$ (strongly if  $p<\infty$) in $L^p(\mu)$. Observe that the function
$f_n=E[f|\sigma(\pi_{F_n})]$ can be identified with a function on the
finite dimensional abstract Wiener space $(F_n,\mu_n,F_n)$, where
$\mu_n=\pi_n\mu$.

Since the translations of $\mu$ with the elements of $H$ induce measures
equivalent to $\mu$, the G\^ateaux  derivative in $H$ direction of the
random variables is a closable operator on $L^p(\mu)$-spaces and  this
closure will be denoted by $\nabla$ cf.,  for example
\cite{F-P},\cite{ASU, ASU-1}. The corresponding Sobolev spaces
(the equivalence classes) of the  real random variables
will be denoted as $\DD_{p,k}$, where $k\in \NN$ is the order of
differentiability and $p>1$ is the order of integrability. If the
random variables are with values in some separable Hilbert space, say
$\Phi$, then we shall define similarly the corresponding Sobolev
spaces and they are denoted as $\DD_{p,k}(\Phi)$, $p>1,\,k\in
\NN$. Since $\nabla:\DD_{p,k}\to\DD_{p,k-1}(H)$ is a continuous and
linear operator its adjoint is a well-defined operator which we
represent by $\delta$.  $\delta$ coincides with the It\^o
integral of the Lebesgue density of the adapted elements of
$\DD_{p,k}(H)$ (cf.\cite{ASU,ASU-1}).

For any $t\geq 0$ and measurable $f:W\to \reals_+$, we note by
$$
P_tf(x)=\int_Wf\left(e^{-t}x+\sqrt{1-e^{-2t}}y\right)\mu(dy)\,,
$$
it is well-known that $(P_t,t\in \reals_+)$ is a hypercontractive
semigroup on $L^p(\mu),p>1$,  which is called the Ornstein-Uhlenbeck
semigroup (cf.\cite{F-P,ASU,ASU-1}). Its infinitesimal generator is denoted
by $-\calL$ and we call $\calL$ the Ornstein-Uhlenbeck operator
(sometimes called the number operator by the physicists). The
norms defined by
\begin{equation}
\label{norm}
\|\phi\|_{p,k}=\|(I+\calL)^{k/2}\phi\|_{L^p(\mu)}
\end{equation}
are equivalent to the norms defined by the iterates of the  Sobolev
derivative $\nabla$. This observation permits us to identify the duals
of the space $\DD_{p,k}(\Phi);p>1,\,k\in\NN$ by $\DD_{q,-k}(\Phi')$,
with $q^{-1}=1-p^{-1}$,
where the latter  space is defined by replacing $k$ in (\ref{norm}) by
$-k$, this gives us the distribution spaces on the Wiener space $W$
(in fact we can take as $k$ any real number). An easy calculation
shows that, formally, $\delta\circ \nabla=\calL$, and this permits us
to extend the  divergence and the derivative  operators to the
distributions as linear,  continuous operators. In fact
$\delta:\DD_{q,k}(H\otimes \Phi)\to \DD_{q,k-1}(\Phi)$ and
$\nabla:\DD_{q,k}(\Phi)\to\DD_{q,k-1}(H\otimes \Phi)$ continuously, for
any $q>1$ and $k\in \reals$, where $H\otimes \Phi$ denotes the
completed Hilbert-Schmidt tensor product (cf., for instance
\cite{Mal,ASU,ASU-1}). Finally, in the case of classical Wiener space,
we denote by $\DD_{p,k}^a(H)$ the subspace defined by 
$$
\DD_{p,k}^a(H)=\{\xi\in\DD_{p,k}(H):\,\dot{\xi}\mbox{ is adapted}\}
$$
for $p\geq 1,\,k\in\R$. 

Let us recall some facts from the convex analysis. Let $K$ be a
Hilbert space, a subset $S$ of $K\times K$ is called cyclically
monotone if  any finite subset
$\{(x_1,y_1),\ldots,(x_N,y_N)\}$ of
$S$ satisfies the following algebraic condition:
$$
\langle y_1,x_2-x_1\rangle+\langle y_2,x_3-x_2\rangle+\cdots+\langle
y_{N-1},x_N-x_{N-1}\rangle+\langle y_N,x_1-x_N\rangle\leq 0\,,
$$
where $\langle\cdot,\cdot\rangle$ denotes the inner product of
$K$. It turns out  that $S$ is
cyclically monotone if and only if
$$
\sum_{i=1}^N(y_i,x_{\sigma(i)}-x_i)\leq 0\,,
$$
for any permutation $\sigma$ of $\{1,\ldots,N\}$ and for any
finite subset $\{(x_i,y_i):\,i=1,\ldots,N\}$ of $S$. Note that $S$
is cyclically monotone if and only if any translate of it is
cyclically monotone.  By a theorem of Rockafellar,  any cyclically
monotone set is contained in the graph of the subdifferential  of
a convex function in the sense of convex analysis (\cite{Roc}) and
even if the function may not be unique its subdifferential  is
unique.

\noindent
Let now  $(W,\mu,H)$ be an abstract Wiener space;  a measurable  function
$f:W\to \reals\cup\{\infty\}$  is called $1$-convex if the map
$$
h\to f(x+h)+\frac{1}{2}|h|_H^2=F(x,h)
$$
is convex on the Cameron-Martin space $H$ with values in
$L^0(\mu)$. Note that this notion is compatible with the
$\mu$-equivalence classes of random variables thanks to the
Cameron-Martin theorem. It is proven in \cite{F-U1} that
this definition  is equivalent  the following condition:
  Let $(\pi_n,n\geq 1)$ be a sequence of regular, finite dimensional,
  orthogonal projections of
  $H$,  increasing to the identity map
  $I_H$. Denote also  by $\pi_n$ its  continuous extension  to $W$ and
  define $\pi_n^\bot=I_W-\pi_n$. For $x\in W$, let $x_n=\pi_nx$ and
  $x_n^\bot=\pi_n^\bot x$.   Then $f$ is $1$-convex if and only if
$$
x_n\to \frac{1}{2}|x_n|_H^2+f(x_n+x_n^\bot)
$$
is  $\pi_n^\bot\mu$-almost surely convex.
\subsection{Preliminaries about the Monge-Kantorovitch measure
  transportation problem}

\begin{definition}
Let $\xi$ and $\eta$ be two probabilities on $(W,\calB(W))$. We say
that a probability $\ga$ on $(W\times W,\calB(W\times W))$ is a
solution of the Monge-Kantorovitch problem associated to the
couple $(\xi,\eta)$ if the first marginal of $\ga$ is $\xi$, the
second one is $\eta$ and if
$$
J(\ga)=\int_{W\times W}|x-y|_H^2d\ga(x,y)=\inf\left\{\int_{W\times
  W}|x-y|_H^2d\beta(x,y):\,\beta\in \Sigma(\xi,\eta)\right\}\,,
$$
where $\Sigma(\xi,\eta)$ denotes the set of all the probability
measures on $W\times W$ whose first and second marginals are
respectively $\xi$ and $\eta$. We shall denote the Wasserstein
distance between $\xi$ and $\eta$, which is  the positive
square-root of  this infimum, with $d_H(\xi,\eta)$.
\end{definition}

\remark
By the weak compacteness  of probability measures on $W\times W$
and the  lower semi-continuity of the strictly convex cost function,
the infimum in the definition is  attained
even if the functional $J$ is identically infinity. In this latter
case we say that the solution is degenerate.

The next result, which is the extension   of
the finite dimensional version of an inequality due to  Talagrand, \cite{Tal},  gives a
sufficient condition for the finiteness of the Wasserstein
distance in the case one of the measures is the Wiener measure $\mu$  and the
second one is absolutely continuous with respect to it. We give a
short proof for the sake of completeness:
\begin{theorem}
\label{ineq-thm}
Let $L\in L\log L(\mu)$ be a positive random variable with
$E[L]=1${\footnote{In the sequel we denote the expectation w.r. to the
  Wiener measure by $E$}} and let  $\nu$ be  the measure $d\nu=Ld\mu$.
We then  have
\begin{equation}
\label{tal-ineq}
d_H^2(\nu,\mu)\leq 2E[L\log L]\,.
\end{equation}
\end{theorem}
\nproof Let us remark first that we can take $W$ as the classical
Wiener space $W=C_0([0,1])$ and, using the stopping techniques of
the martingale theory, we may assume that $L$ is upper and lower
bounded almost surely. Then a classical result of the It\^o
calculus implies that $L$ can be represented as an exponential
martingale
$$
L_t=\exp\left\{-\int_0^t\dot{u}_\tau
  dW_\tau-\frac{1}{2}\int_0^t|\dot{u}_\tau|^2d\tau\right\}\,,
$$
with $L=L_1$, where $(\dot{u}_t,t\in [0,1])$ is a measurable  process
adapted to the filtration of the canonical  Wiener process $(t,x)\to
W_t(x)=x(t)$. Let us define  $u:W\to H$  as
$u(t,x)=\int_0^t\dot{u}_\tau(x) d\tau$ and $U:W\to W$ as
$U(x)=x+u(x)$. The Girsanov theorem implies that $x\to U(x)$ is a
Browian motion under $\nu$, hence the image of the measure $\nu$ under
the map $U\times I_W:W\to W\times W$ denoted by $\beta=(U\times I)\nu$
belongs to
$\Sigma(\mu,\nu)$. Let $\ga$ be any optimal measure, then
\beaa
J(\ga)&=&d_H^2(\nu,\mu)\leq \int_{W\times W}|x-y|_H^2d\beta(x,y)\\
&=&E[|u|_H^2L]\\
&=&2E[L\log L]\,, \eeaa where the last equality follows also from
the Girsanov theorem and the It\^o stochastic calculus.
\nqed

\noindent
The next two theorems, which explain the existence and several
properties of the solutions of Monge-Kantorovitch problem and the
transport maps  have been proven in \cite{F-U2}.
\begin{theorem}[General case]
\label{monge-general}
Suppose that $\rho$ and $\nu$ are two probability measures on
$W$ such  that
$$
d_H(\rho,\nu)<\infty\,.
$$
Let $(\pi_n,n\geq 1)$ be a total increasing  sequence of regular
projections (of $H$, converging to the identity map of $H$).
Suppose  that, for any $n\geq 1$, the regular
conditional probabilities $\rho(\cdot\,|\pi_n^\bot=x^\bot)$ vanish
$\pi_n^\bot\rho$-almost surely on
the subsets of  $(\pi_n^\bot)^{-1}(W)$ with Hausdorff dimension
$n-1$. Then there exists a
unique solution of the   Monge-Kantorovitch problem, denoted by $\ga\in
\Sigma(\rho,\nu)$ and  $\ga$ is supported by the graph of a Borel
map $T$ which is the solution of the
Monge problem.  $T:W\to W$ is  of the form $T=I_W+\xi$ , where $\xi\in
H$ almost surely. Besides  we have
\beaa
d_H^2(\rho,\nu)&=&\int_{W\times W}|T(x)-x|_H^2d\ga(x,y)\\
&=&\int_{W}|T(x)-x|_H^2d\rho(x)\,,
\eeaa
and
for $\pi_n^\bot\rho$-almost almost all $x_n^\bot$, the map $u\to
\xi(u+x_n^\bot)$ is cyclically monotone on
$(\pi_n^\bot)^{-1}\{x_n^\bot\}$, in the sense that
$$
\sum_{i=1}^N\left(\xi(x_n^\bot+u_i),u_{i+1}-u_i\right)_H\leq 0
$$
$\pi_n^\bot\rho$-almost surely, for any cyclic sequence
$\{u_1,\ldots,u_N,u_{N+1}=u_1\}$ from $\pi_n(W)$. Finally, if, for any $n\geq
1$, $\pi_n^\bot\nu$-almost surely,  $\nu(\cdot\,|\pi_n^\bot=y^\bot)$
 also vanishes on the $n-1$-Hausdorff dimensional  subsets
 of $(\pi_n^\bot)^{-1}(W)$, then $T$ is invertible, i.e, there exists
 $S:W\to W$ of the form $S=I_W+\eta$ such that  $\eta\in H$ satisfies
 a similar  cyclic monotononicity property as $\xi$ and that
\beaa
1&=&\ga\left\{(x,y)\in W\times W: T\circ S(y)=y\right\}\\
&=&\ga\left\{(x,y)\in W\times W: S\circ T(x)=x\right\}\,.
\eeaa
In particular we have
\beaa
d_H^2(\rho,\nu)&=&\int_{W\times W}|S(y)-y|_H^2d\ga(x,y)\\
&=&\int_{W}|S(y)-y|_H^2d\nu(y)\,.
\eeaa
\end{theorem}
\begin{remarkk}
In particular, for all  the measures $\rho$ which are  absolutely
continuous with respect to the  Wiener measure $\mu$,  the second
hypothesis is satisfied, i.e., the measure
$\rho(\cdot\,|\pi_n^\bot=x_n^\bot)$ vanishes on the sets of Hausdorff
dimension $n-1$.
\end{remarkk}

\noindent
The case where one of the measures is the Wiener measure  and the other is
absolutely continuous with respect to $\mu$ is the most important one
for the applications. Consequently we give the related results
separately in the following theorem where the tools of the Malliavin
calculus give more information about the maps $\xi$ and $\eta$ of
Theorem \ref{monge-general}:
\begin{theorem}[Gaussian case]
\label{gaussian-case}
Let $\nu$ be the measure $d\nu=Ld\mu$, where $L$ is a positive random variable,
with $E[L]=1$. Assume that $d_H(\mu,\nu)<\infty$ (for instance
$L\in L\log L$). Then there exists a  $1$-convex function $\phi\in
\DD_{2,1}$, unique up to a constant,  such that  the map
$T=I_W+\nabla \phi$ is the unique solution of the original problem of
Monge. Moreover, its graph supports  the  unique
solution of the
Monge-Kantorovitch problem $\ga$. Consequently
$$
(I_W\times T)\mu=\ga
$$
In particular  $T$ maps $\mu$ to $\nu$ and  $T$ is almost surely
invertible, i.e., there exists some $T^{-1}$ such that $T^{-1}\nu=\mu$
and that
\beaa
1&=&\mu\left\{x:\,T^{-1}\circ T(x)=x\right\}\\
&=&\nu\left\{y\in W:\,T\circ T^{-1}(y)=y\right\}\,.
\eeaa
\end{theorem}

\begin{remarkk}
\label{nu-closed}
Assume that the operator $\nabla$ is closable with respect to $\nu$,
then we have $\eta=\nabla\psi$. In particular, if  $\nu$ and $\mu$ are
equivalent, then we have
$$
T^{-1}=I_W+\nabla\psi\,,
$$
where is $\psi$ is  a $1$-convex function. $\psi$ is called the dual
potential of the MKP$(\mu,\nu)$ and we have the following relations:
$$
\phi(x)+\psi(y)+\frac{1}{2}|x-y|_H^2\geq 0\,,
$$
for any $x,y\in W$, and
$$
\phi(x)+\psi(y)+\frac{1}{2}|x-y|_H^2=0
$$
$\ga$-almost surely.
\end{remarkk}
\label{stability}
\begin{remarkk}
Let $(e_n,n\in \NN)$ be a complete,  orthonormal in $H$, denote by
$V_n$ the sigma algebra generated by $\{\delta e_1,\ldots,\delta
e_n\}$ and let  $L_n=E[L|V_n]$. If $\phi_n\in \DD_{2,1}$ is the
function constructed in Theorem \ref{gaussian-case}, corresponding to
$L_n$, then, using the inequality (\ref{tal-ineq})  we can prove that
the sequence $(\phi_n,n\in \NN)$ converges to $\phi$ in $\DD_{2,1}$.
\end{remarkk}

\section{Characterization of the invertible shifts}
Let us begin with some results of general interest. 
Let us first define:
\begin{definition}
A measurable map $T:W\to W$ is called ($\mu$-) almost surely right invertible
if there exists a measurable  map $S:W\to W$ such that $S\mu\ll \mu$
and $T\circ S=I_W$ $Óu$-a.s. Similarly, we say that it is left
invertible, if $T\mu\ll \mu$ and if there exists a measurable map
$S:W\to W$ such that $S\circ T=I_W$$\mu$-a.s.
\end{definition}

\noindent
The following
proposition some parts of which  are  proven in \cite{INV}, shows that, whenever an adapted
shift has a left inverse almost surely, then it is almost surely
invertible and  its inverse is also an adapted perturbation of
identity  and it relates this concept to the existence and uniqueness
of strong solutions of
stochastic differential equations. The a
\begin{proposition}
\label{inverse-1}
Assume $A=I_W+a$, $a\in L^2(\mu,H)$, $\dot{a}$ is adapted,
$E[\rho(-\delta a)]=1$. Suppose that there exists a map  $B:W\to W$ such
that $B\circ A=I_W$ a.s. Then the following assertions are true:
\begin{itemize}
\item[(i)] $B\mu$ is equivalent to $\mu$ and $A\circ B=I_W$ a.s., i.e.,
  $B$ is also a right inverse.
\item[(ii)] $B=I_W+b$, $b:W\to H$, $\dot{b}$ is also adapted.
\item[(iii)] $(t,w)\to B_t(w)$ is the strong solution of
\begin{eqnarray}
dB_t&=&-\dot{a}_t\circ Bdt+dW_t\label{sde}\\
B_0&=&0\nonumber\,.
\end{eqnarray}
\item[(iv)] We have 
\begin{eqnarray}
\dot{a}_t+\dot{b}_t\circ A&=&0\label{C1}\\
\dot{b}_t+\dot{a}_t\circ B&=&0\label{C2}
\end{eqnarray}
$dt\times d\mu$-a.s.
\item[(v)] In particular either  the property  $A\mu\sim \mu$ and the relation 
 \ref{C1}  together or  $B\mu\sim \mu$ and the relation 
 \ref{C2}  together imply  that $B\circ A=A\circ B=I_W$ a.s.
\end{itemize}
\end{proposition}
\nproof For any $f\in C_b(W)$, it follows  from the Girsanov
theorem \beaa
E[f\circ B]&=&E[f\circ B\circ A\,\rho(-\delta a)]\\
&=&E[f\,\rho(-\delta a)]\,,
\eeaa
hence $B\mu$ is equivalent to
$\mu$ and the corresponding Radon-Nikodym density is $\rho(-\delta
a)$. Let
$$
D=\{w\in W:\,B\circ A(w)=w\}\,.
$$
Since $D\subset A^{-1}(A(D))$ and by  the hypothesis $\mu(D)=1$ we
get
$$
E[1_{A(D)}\circ A]=1\,.
$$
Since $A\mu$ is equivalent to $\mu$ we have also $\mu(A(D))=1$. If
$w\in A(D)$, then $w=A(d)$, for some $d\in D$, hence $A\circ
B(w)=A\circ B\circ A(d)=A(d)=w$, consequently $A\circ B=I_W$
$\mu$-almost surely and $B$ is the two-sided inverse of $A$.
Evidently, together with the absolute continuity of $B\mu$, this
implies that $B$ is of the form $B=I_W+b$, with $b:W\to H$.
Moreover, $\dot{a}=-\dot{b}\circ A$, hence the right hand side is
adapted. We can assume that all these processes are
uni-dimensional (otherwise we proceed component wise). Let
$\dot{b}^n=\max(-n,\min(\dot{b},n))$. Then $\dot{b}^n\circ A$ is
adapted. Let $H\in L^2(dt\times d\mu)$ be an adapted process.
Using the Girsanov theorem: 
\beaa
 E\left[\rho(-\delta
a)\int_0^1\dot{b}^n_s\circ A\,H_s\circ A
ds\right]&=&E\left[\int_0^1\dot{b}^n_sH_sds\right]\\
&=&E\left[\int_0^1E[\dot{b}^n_s|\calF_s]H_sds\right]\\
&=&E\left[\rho(-\delta a)\int_0^1E[\dot{b}^n_s|\calF_s]\circ
A\,H_s\circ A ds\right]\,. \eeaa Consequently
$$
E[\dot{b}^n_s|\calF_s]\circ A=\dot{b}^n_s\circ A\,,
$$
almost surely. Since $A\mu$ is equivalent to $\mu$, it follows
that
$$
E[\dot{b}^n_s|\calF_s]=\dot{b}^n_s
$$
almost surely, hence $\dot{b}^n$ and consequently  $\dot{b}$ are
adapted. It is now clear that $(B(t),\,t\in [0,1])$ is a strong
solution of (\ref{sde}). The uniqueness follows from the fact
that, any strong solution of (\ref{sde}) would be a right  inverse
to $A$, since $A$ is invertible, then this solution is equal to
$B$. 

The proof of (v) is quite similar to that of the first part: let
$D=\{w\in W: A\circ B(w)=w\}$, then $\mu(B^{-1}(B(D))=1$, hence $B\circ
A=I_W$$\mu$-a.s. Moreover $B$ can be written as $B=I_W+b$, with
$\dot{a}=-\dot{b}\circ A$, proceeding as above, we show that $\dot{b}$ is adapted and the rest of the
proof follows.
\nqed

 The invertibility of $A$ is characterized in terms of the
 corresponding Wick exponentials  as below:
\begin{theorem}
\label{inverse-thm}
Let  $A=I_W+a$, $a\in L^0_a(\mu,H)$.
Assume  that $E[\rho(-\delta a)]=1$\footnote{Here we denote by $\delta
  a$ the stochastic integral of the adapted process $\dot{a}$ in
 $L^0(\mu)$} and that
$$
\frac{dA\mu}{d\mu}\circ A\,\rho(-\delta a)=1
$$
almost surely. Then $A$ is (almost surely) invertible.
\end{theorem}
\nproof Since $E[\rho(-\delta a)]=1$, $A\mu$ is equivalent to
$\mu$, hence the corresponding Radon-Nikodym derivative can be
expressed as an exponential martingale:
$$
l=\frac{dA\mu}{d\mu}=\exp\left(-\delta b-\frac{1}{2}|b|_H^2\right)\,,
$$
where $b(t,w)=\int_0^t\dot{b}_s(w)ds$, with $\dot{b}$ adapted,
$\int_0^1|\dot{b}_s|^2ds<\infty$ almost surely and $\delta b$ is
defined in $L^0(\mu)$. The hypothesis implies that
\begin{equation}
\label{H1}
\delta(a+b\circ A)+\frac{1}{2}|a+b\circ A|_H^2=0
\end{equation}
almost surely. Define the local martingale $(M_t)$ as
$$
M_t=\exp\left(-\int_0^t(\dot{a}_s+\dot{b}_s\circ
  A)dW_s-\frac{1}{2}\int_0^t|\dot{a}_s+\dot{b}_s\circ A|^2ds\right)\,.
$$
The relation (\ref{H1}) implies in fact that $(M_t)$ is a uniformly
integrable martingale with its final value (at $t=1$)
$M_1=1$. Consequently $M_t=1$ almost surely for any $t\in [0,1]$ and
this implies that
$$
\dot{a}_s+\dot{b}_s\circ A=0
$$
$ds\times d\mu$-almost surely. Hence $(I_W+b)\circ A=I_W$ almost
surely and the proof is fully completed thanks to Proposition
\ref{inverse-1}. \nqed

\begin{proposition}
Assume that $(A_n, n\geq 1)$ is a sequence of mappings of the form
$A_n=I_W+a_n$, with $a_n:W\to H$, $\dot{a}_n$ is adapted for any $n$
and $(a_n, n\geq 1)$ converges to some $a$ in $L^0(\mu,H)$ such that
$E[\rho(-\delta a)]=1$. Suppose
that, for any $n\geq 1$,  $E[\rho(-\delta a_n)]=1$ and $A_n$ is invertible. If
$$
\lim_{n\to \infty}\frac{dA_n\mu}{d\mu}=l
$$
exists in the norm topology of $L^1(\mu)$, then $A=I_W+a$ is also
invertible.
\end{proposition}
\nproof Let us denote by $l_n$ the Radon-Nikodym derivative of
$A_n\mu$ with respect to $\mu$. The hypothesis implies that $(l_n,
n\geq 1)$ is uniformly integrable. Since $(a_n, n\geq 1)$
converges in probability, the uniform integrability, combined with
the Lusin theorem implies that $(l_n\circ A_n,n\geq 1)$ converges
in probability to $l\circ A$. Since $(\rho(-\delta a_n),n\geq 1)$
converges to $\rho(-\delta a)$ in probability and since, by the
invertibility of $A_n$, we have
$$
l_n\circ A_n\,\rho(-\delta a_n)=1
$$
almost surely for any $n\geq 1$, we have also
$$
l\circ A\,\rho(-\delta a)=1
$$
almost surely. The conclusion follows then from Theorem
\ref{inverse-thm}. \nqed

The following lemma gives an important information about the
Radon-Nikodym  density of the measure $A\mu$ with respect to $\mu$:
\begin{lemma}
\label{R-N-lemma}
Assume that $A=I_W+a$ with $a\in L^0(\mu, H)$ with $\dot{a}$
adapted. Then
$$
\frac{dA\mu}{d\mu}\circ A\,\,E[\rho(-\delta a)|A]\leq 1
$$
almost surely. If we have also
$E[\rho(-\delta a)]=1$, then the above inequality becomes an equality:
$$
\frac{dA\mu}{d\mu}\circ A\,\,E[\rho(-\delta a)|A]=1
$$
almost surely.
\end{lemma}
\nproof For any positive function  $f\in C_b(W)$, using the
Girsanov theorem and the Fatou Lemma, we have \beaa
 E[f\circ A]&=&E\left[f\,\frac{dA\mu}{d\mu}\right]\\
&\geq&E\left[f\circ A\,\frac{dA\mu}{d\mu}\circ A\,\rho(-\delta a)\right]\\
&=&E\left[f\circ A\,\frac{dA\mu}{d\mu}\circ A\,E[\rho(-\delta
a)|A]\right]\,, \eeaa which proves the first part of the lemma.
For the second part, due to the integrability hypothesis, we can
replace the inequality above by the equality
and the proof follows. \nqed

\section{Properties of non-invertible adapted perturbation of identity}
\noindent
In this section we study the following concept:
\begin{definition}
A positive random variable whose expectation is equal to one with
respect to Wiener measure is said to be representable with a mapping
$U:W\to W$ if 
$$
\frac{dU\mu}{d\mu}=L\,.
$$
\end{definition}
We begin with the following 
\begin{proposition}
Assume that $L=\rho(-\delta v)$, where $v\in L_a^0(\mu,H)$, i.e.,
$\dot{v}$ is adapted and $\int_0^1|\dot{v}_s|^2ds<\infty$ a.s.
Then there exists $U=I_W+u$, with $u:W\to H$ adapted such that
$U\mu=L\mu$ and $E[\rho(-\delta u)]=1$ if and only if the
following condition is satisfied:
\begin{eqnarray}
\label{MA-1}
1&=&L_t\circ U \,\,E\left[\rho(-\delta u^t)|\calU_t\right]\\
&=&L_t\circ U \,\,E\left[\rho(-\delta u)|\calU_t\right]
\end{eqnarray}
almost surely for any $t\in [0,1]$, where $u^t$ is defined as
$u^t(\tau)=\int_0^{t\wedge \tau}\dot{u}_sds$  and  $\calU_t$ is
the sigma algebra generated by $(w(\tau)+u(\tau),\,\tau\leq t)$.
\end{proposition}
\nproof Let $U_t$ be defined as $I_W+u^t$, then for any $f\in
C_b(W)$ which is $\calF_t$-measurable, we have \beaa
E[f\circ U_t\,L_t\circ U_t\,\rho(-\delta u^t)]&=&E[f\,L_t]\\
&=&E[f\circ U_t]\,. \eeaa Since, for any  $\calF_t$-measurable
function $G$, $G\circ U_t$ is $\calU_t$ measurable, we get
$$
L_t\circ U_t \,\,E\left[\rho(-\delta u^t)|\calU_t\right]=1\,.
$$
Conversely, it follows from the relation (\ref{MA-1}) and from the
Girsanov theorem that
$$
E[f\circ U]=E[f\circ U\,L\circ U\,\rho(-\delta u)]=E[f\,L]\,,
$$
a similar relation holds when we replace $U$ by $U_t$. \nqed

\noindent
Let us calculate $E[\rho(-\delta
u^t)|\calU_t]=E[\rho(-\delta u)|\calU_t]$ in terms of the
innovation process associated to $U$. Recall that the term
innovation, which originates from the filtering theory is defined
as (cf.\cite{FKK} and  \cite{BOOK})
$$
Z_t=U_t-\int_0^t E[\dot{u}_s|\calU_s]ds
$$
and it is a $\mu$-Brownian motion with respect to the filtration
$(\calU_t,t\in [0,1])$. A similar proof as the one in \cite{FKK} shows
that any martingale with respect to the filtration of $U$ can be
represented as a stochastic integral with respect to $Z$. Hence, by
the positivity assumption, $E[\rho(-\delta u)|\calU_t]$ can be written
as an exponential martingale
$$
E[\rho(-\delta
u)|\calU_t]=\exp\left(-\int_0^t(\dot{\xi}_s,dZ_s)-\frac{1}{2}\int_0^t|\dot{\xi}_s|^2ds\right)\,.
$$
Below we give a more detailed result:
\begin{proposition}
\label{iyi-prop}
We have the following explicit result
\begin{equation}
\label{cond-exp} E[\rho(-\delta
u)|\calU]=\exp\left(-\int_0^1(E[\dot{u}_s|\calU_s],dZ_s)-\frac{1}{2}\int_0^1|E[\dot{u}_s|\calU_s]|^2ds\right)\,,
\end{equation}
hence 
\begin{equation}
\label{cond-exp1} E[\rho(-\delta
u)|\calU_t]=\exp\left(-\int_0^t(E[\dot{u}_s|\calU_s],dZ_s)-\frac{1}{2}\int_0^t|E[\dot{u}_s|\calU_s]|^2ds\right)\,,
\end{equation}
almost surely.
\end{proposition}
\nproof The proof follows from the double utilization of the
Girsanov theorem. Let us denote by $l_t$ the Girsanov exponential 
$$
l_t=\exp\left(-\int_0^t(E[\dot{u}_s|\calU_s],dZ_s)-\frac{1}{2}\int_0^t|E[\dot{u}_s|\calU_s]|^2ds\right)\,.
$$
 On the first hand, we have, for any $f\in C_b(W)$,
$$
E[f\circ U\rho(-\delta u)]=E[f]\,,
$$
and on the other hand, applying the Girsanov theorem to  the decomposition
$$
U_t=Z_t+\int_0^t E[\dot{u}_s|\calU_s]ds\,,
$$
we get 
$$
E[f\circ U\,l_1]\leq E[f]=E[f\circ U\, \rho(-\delta u)]
$$
for any positive, measurable $f$ on $W$. Taking $f$ to be $\calF_t$
measurable, we conclude that 
$$
l_t\leq E[\rho(-\delta u)|\calU_t]
$$
a.s. for any $t\in[0,1]$. Consequently $(l_t,t\in[0,1])$ is a
uniformly integrable martingale and in particular $E[l_1]=1$. Hence we
have 
$$
E\left[f\circ U\,l_1\right]=E[f]=E[f\circ
U\,\rho(-\delta u)]\,,
$$
for any $f\in C_b(W)$ which implies that $l_1=E[\rho(-\delta
u)|\calU]$  and the proof of (\ref{cond-exp}) follows.
The relation (\ref{cond-exp1}) is obvious  since
$\calU_t\subset\calF_t$.

 \nqed

\begin{theorem}
\label{thm-1}
A necessary and sufficient condition for the relation (\ref{MA-1}),
that is to say for the representability of $L=\rho(-\delta v)$ by $U=I_W+u$
is that
$$
E[\dot{u}_t|\calU_t]=-\dot{v}_t\circ U
$$
$dt\times d\mu$-almost surely.
\end{theorem}
\nproof We have
$$
L_t\circ U=\exp\left(-\delta v_t\circ U_t-\frac{1}{2}|v_t\circ
  U_t|_H^2\right)\,.
$$
Moreover using the identity
$$
\delta v_t\circ U_t=\int_0^t(\dot{v}_s\circ
U,dW_s)+\int_0^t(\dot{v}_s\circ U,\dot{u}_s)ds\,,
$$
we get
$$
L_t\circ
U=\exp\left[-\int_0^t\left(\dot{v}_s\circ
  U,dW_s+\dot{u}_sds+\frac{1}{2}\dot{v}_s\circ U\,ds\right)\right]\,.
$$
Substituting all these relations in  (\ref{MA-1}) and using the
representation (\ref{cond-exp}), we obtain
\beaa
1&=&L_t\circ U\,\,E[\rho(-\delta u)|\calU_t]\\
&=&\exp\left[-\int_0^t\left(\dot{v}_s\circ
  U,dW_s+\dot{u}_sds+\frac{1}{2}\dot{v}_s\circ U\,ds\right)\right]\\
&&\exp\left(-\int_0^t(E[\dot{u}_s|\calU_s],dZ_s)-
\frac{1}{2}\int_0^t|E[\dot{u}_s|\calU_s]|^2ds\right)\,.
\eeaa
But
$$
\int_0^t(E[\dot{u}_s|\calU_s],dZ_s)=\int_0^t(E[\dot{u}_s|\calU_s],dW_s+(\dot{u}_s-
E[\dot{u}_s|\calU_s])ds)\,.
$$
Consequently we get
$$
\int_0^t(\dot{v}_s\circ U+E[\dot{u}_s|\calU_s],dW_s)=0\,,
$$
almost surely for any $t\in[0,1]$ and this implies that
$$
E[\dot{u}_s|\calU_s]=-\dot{v}_s\circ U
$$
$ds\times d\mu$-almost surely. The sufficiency is obvious.
 \nqed

\begin{corollary}
\label{d=1-cor}
A  necessary and sufficient condition for the
relation (\ref{MA-1}) is that
that
$$
V\circ U=Z\,,
$$
in other words
$$
U_t=Z_t-\int_0^t\dot{v}_s\circ U\,ds
$$
almost surely, where $Z$ is the innovation process associated to $U$.
\end{corollary}
\nproof 
The condition in Theorem \ref{thm-1} reads as
\begin{equation}
\label{nice-cond}
\dot{v}_t\circ U+E[\dot{u}_t|\calU_t]=0
\end{equation}
almost surely. Hence
\beaa
(V\circ U)(t)&=&U(t)+(v\circ U)(t)\\
&=&Z(t)+\int_0^tE[\dot{u}_s|\calU_s]ds+\int_0^t\dot{v}_s\circ U\,ds\\
&=&Z_t\,, \eeaa by the relation (\ref{nice-cond}). \nqed

\begin{corollary}
\label{mloc-cor}
Suppose that the innovation process $Z$ is an
$(\calF_t,t\in[0,1])$-local martingale, then $U$ is almost surely
invertible and its inverse is $V$.
\end{corollary}
\nproof We have
$$
U_t=W_t+\int_0^t\dot{u}_sds=Z_t+\int_0^tE[\dot{u}_s|\calU_s]ds\,,
$$
hence
$(W_t-Z_t,t\in [0,1])$ is a continuous local  martingale of finite
variation. This implies that $Z$ and $W$ are equal hence
$$
\dot{u}_t=E[\dot{u}_t|\calU_t]\,,
$$
$dt\times d\mu$-almost surely. From Theorem \ref{thm-1}, it follows
that $u+v\circ U=0$ almost surely, i.e., $V\circ U=I_W$ almost
surely. It follows from Proposition \ref{inverse-1} that
$$
U\circ V=I_W
$$
also $\mu$-almost surely. \nqed

\noindent
 We can give a complete characterization of the
representable random variables as follows:
\begin{theorem}
Assume that $L=\rho(-\delta v)$, $V=I_W+v$, $v\in L^0_a(\mu,H)$. Assume that $U=I_W+u$
is also an adapted perturbation of identity with $E[\rho(-\delta
u)]=1$. Assume that $V\circ U=B$ is a Brownian motion with respect
to its own filtration. We have $U\mu=L\cdot\mu$ if and only if $B$
is a local martingale with respect to the filtration generated by
$U$ and in this case $B$ is equal to the innovation associated to $U$.
\end{theorem}
\nproof 
The necessity has already been proven, for the
sufficiency, note that,  we have $U=B-v\circ U$. On the other
hand we can always represent $U$ by its innovation process as
$$
U_t=Z_t+\int_0^t E[\dot{u}_s|\calU_s]ds=B_t-\int_0^t\dot{v}_s\circ Uds
$$
where $Z$ is the innovation process associated to $U$, which is  a Brownian motion with respect to $(\calU_t,\,t\in
[0,1])$. Consequently
$$
-\dot{v}_s\circ U=E[\dot{u}_s|\calU_s]\,,
$$
$ds\times d\mu$-almost surely and the proof follows from Theorem
\ref{thm-1}. \nqed

\section{Relations with entropy}
Assume that $u\in\DD_{2,0}^a(H)$ with $E[\rho(-\delta u)]=1$ and
let $L\in \L\log\L(\mu)$ be the Radon-Nikodym density of
$U\mu=(I_W+u)\mu$ with respect to $\mu$. Let us represent $L$ as
$\rho(-\delta v)$. Denote $E[\rho(-\delta u)|\calU]$ by
$\hat{\rho}$. Then, due to the Girsanov theorem, we have \beaa
E[\hat{\rho}\log\hat{\rho}]&=&\half E\left[\hat{\rho}|v\circ U|_H^2\right]\\
&=&\frac{1}{2}E\left[\rho(-\delta u)|v\circ U|_H^2\right]\\
&=&\frac{1}{2}E[|v|_H^2]\,.
\eeaa
In particular, the Jensen inequality implies that
\beaa
E[|v|_H^2]&\leq &2 E[\rho(-\delta u)\log \rho(-\delta u)]\\
&=&E[\rho(-\delta u)|u|_H^2]\,.
\eeaa

\begin{proposition}
 Let $P_\eps$ denote the
Ornstein-Uhlenbeck semigroup and denote by $v_\eps$ the regularization
$P_\eps v$ and denote by $u_\eps$ the $H$-valued mapping which is
defined as $I_W+u_\eps=(I_W+v_\eps)^{-1}$ whose  existence
 follows from \cite{INV}. The set $(u_\eps,\eps>0)$ has a unique
 weak  accumulation point  $\tilde{u}\in \DD_{2,0}(H)$. If the relation
(\ref{MA-1}) holds  then $\tilde{u}$ satisfies the following relation:
$$
\frac{d}{ds}\tilde{u}(s)\circ Z=-E[\dot{v}_s\circ
U|\calZ_s]=E[\dot{u}_s|\calZ_s]
$$
$ds\times d\mu$-almost surely, where $\calZ$ denotes the sigma algebra
generated by the innovation $Z$ associated to $U$.
\end{proposition}
\nproof
 From \cite{INV},  $V_\eps=I_W+v_\eps$ is almost surely invertible and its
 inverse
 can be written as $U_\eps=I_W+u_\eps$. Moreover $u_\eps=-v_\eps\circ
 U_\eps$. Hence $(u_\eps,\eps>0)$ is bounded in
 $L^2(\mu,H)$. Consequently, there exists a subnet which converges
 weakly to some $\tilde{u}$. Let $\xi$ be an $H$-valued,  bounded
 continuous function on $W$. Denoting by $\lip\cdot,\cdot\rip$ the duality
 bracket of $L^2(\mu,H)$, we get
\beaa
\lip u_\eps,\xi\rip&=&\lip u_\eps\circ V_\eps,\xi\circ
V_\eps\,\rho(-\delta v_\eps)\rip\\
&=&-\lip v_\eps, \xi\circ
V_\eps\,\rho(-\delta v_\eps)\rip\\
&\to&-\lip v, \xi\circ V\,\rho(-\delta v)\rip\,.
\eeaa
Hence
$$
\lip \tilde{u},\xi\rip=-\lip v, \xi\circ V\,\rho(-\delta v)\rip\,.
$$
Consequently $\tilde{u}$ is unique, i.e., the net $(u_\eps,\eps>0)$
has only one accumulation point in the weak topology of
$\DD_{2,0}(H)=L^2(\mu,H)$.
From the last  hypothesis
$$
\frac{dU\mu}{d\mu}=\rho(-\delta v)\,.
$$
Hence
\beaa
\lip \tilde{u},\xi\rip&=&-\lip v, \xi\circ V\,\rho(-\delta v)\rip\\
&=&-\lip v\circ U,\xi\circ V\circ U\rip\\
&=&-\lip v\circ U,\xi\circ Z\rip\\
&=&-E\int_0^1E[\dot{v}_s\circ U|\calZ_s]\,\dot{\xi}_s\circ Z\,ds\,.
\eeaa
Since $Z$ is a Brownian motion, we also have
$$
\lip \tilde{u},\xi\rip=\lip \tilde{u}\circ Z,\xi\circ Z\rip\,,
$$
hence the proof is completed. \nqed
\begin{remarkk}
We draw the attention of the reader to the fact that in general the
weak convergence does not imply the strong convergence. The
situation illustrated above is a typical example for this; in fact if
there were also a strong convergence, then $I+v$ would have been
invertible and we would have $I_W+\tilde{u}=I_W+u=(I_W+v)^{-1}$
(cf. \cite{INV}).
\end{remarkk}
\begin{remarkk}
Similarly, suppose that $v$ is bounded and that
\begin{equation}
\label{hypo}
E[|\tilde{u}|_H^2]=2E[L\log L]\,.
\end{equation}
Then $V=I_W+v$ is invertible and and its inverse is $U=I_W+u$ with
$u=\tilde{u}$. In fact this follows from the hypothesis (\ref{hypo}),
which implies that
\beaa
\lim_{\eps\to 0}E[|u_\eps|_H^2]&=&\lim_{\eps\to 0}E[|v_\eps|_H^2L_\eps]\\
&=&E[|v|_H^2L]\\
&=&2E[L\log L]\\
&=&E[|\tilde{u}|_H^2]\,.
\eeaa
Since $\DD_{2,0}(H)$ is a Hilbert
space, the convergence of the norms implies that  $\lim_{\eps\to
  0}u_\eps=\tilde{u}$ in the norm
topology of $\DD_{2,0}(H)$. Therefore  $V$ is invertible as proven in
\cite{INV}. Consequently, in the case where the mapping $V$ is not
invertible, this equality can not take place.
\end{remarkk}

 \noindent
 The remark above suggests the following claim:
\begin{theorem}
 \label{entropy-thm}
 Assume that $u\in\DD_{2,0}^a(H)$,
$E[\rho(-\delta u)]=1$ and
$$
\frac{dU\mu}{d\mu}=\rho(-\delta v)=L\,,
$$
such that $v\in L_a^0(\mu,H)$. $U=I_W+u$ is then almost surely
invertible with its inverse  $V=I_W+v$ if and only if
$$
2E[L\log L]=E[|u|_H^2]\,.
$$
In other words, $U$ is invertible if and only if
$$
H(U\mu|\mu)=\half \|u\|^2_{\DD_{2,0}(H)}\,,
$$
where $H(U\mu|\mu)$ denotes the entropy of $U\mu$ with respect to
$\mu$.
\end{theorem}
\nproof
 Since $U$ represents $Ld\mu$, we have
$E[\dot{u}_s|\calU_s]+\dot{v}_s\circ U=0$ $ds\times d\mu$-almost
surely. Hence, from the Jensen inequality $E[|v\circ U|_H^2]\leq
E[|u|_H^2]$. Moreover the Girsanov theorem gives
$$
2E[L\log L]=E[|v|_H^2 L]=E[|v\circ
U|_H^2]=E[\int_0^1|E[\dot{u}_s|\calU_s]|^2ds]\,.
$$
Hence the hypothesis implies that
$$
E[|u|_H^2]=E[\int_0^1|E[\dot{u}_s|\calU_s]|^2ds]\,.
$$
From which we deduce that $\dot{u}_s=E[\dot{u}_s|\calU_s]$
$ds\times d\mu$-almost surely. Finally we get
$\dot{u}_s+\dot{v}_s\circ U=0$ $ds\times d\mu$, which is a
necessary and sufficient condition for the claim. The necessity is
obvious.

 \nqed
\begin{remarkk}
\label{rem}
This theorem says that $U$ is invertible if and only if the
``kinetic energy'' of $U$ is equal to the entropy of the measure
that it induces. Moreover $U$ is non-invertible if and only if we
have
$$
H(U\mu|\mu)<\half \|u\|^2_{\DD_{2,0}(H)}\,.
$$
\end{remarkk}

\noindent
The above relation between the entropy and the (kinetic) energy can be
generalized to the maps $I_W+u$, where $u\in L^0(\mu,H)$ which do not
fulfill necessarily the integrability condition $E[\rho(-\delta u)]=1$
as follows:
\begin{theorem}
\label{gen-ineq}
Assume that $u\in L^2_a(\mu,H)$, let $U=I_W+u$ and define $L$ as to be 
$$
L=\frac{dU\mu}{d\mu}\,.
$$
We then  have 
$$
H(U\mu|\mu)=E[L\log L]\leq \half E[|u|_H^2]\,.
$$
\end{theorem}
\nproof
If $|u|_H\in L^\infty(\mu)$, the claim is obvious from above. For the
general case, let $(T_n,n\geq 1)$ be a sequence of stopping times
increasing to infinity such 
that $|u^n|_H$ is bounded, where
$u^n(t)=\int_0^t1_{[0,T_n]}(s)\dot{u}_sds$. Denote by $L_n$ the
Radon-Nikodym derivative of $(I_W+u^n)\mu$ w.r.to $\mu$. From Remark
\cite{rem}, it follows that the sequence $(L_n,n\geq 1)$ is uniformly
integrable, hence it converges to $L$ in the weak topology of
$L^1(\mu)$. From the lower semi-continuity of the entropy w.r. to this
topology, we get
$$
E[L\log L]\leq \lim\inf_n E[L_n\log L_n]\leq \lim \half
E[|u^n|_H^2]=\half E[|u|_H^2]\,.
$$
\nqed

\section{Relations with the innovation conjecture of the filtering}

\noindent 
Let us briefly explain the question (cf. \cite{Zak}, \cite{A-M,FKK} for further
details): Assume that we are given a process of the form
$$
y_t(w,\beta)=W_t(w)+\int_0^t h_s(w,\beta)ds\,,
$$
called the observation, where  $\beta$ is independent of the Wiener path $w$, $s\to
h_s(w,\beta)\in L^2([0,1],ds)$ almost surely and adapted to some
filtration in which the filtration of $(W_t)$ can be injected. The
question is whether the filtration of $y=(y_t,t\in [0,1])$ is equal to
the filtration of the innovation process defined as before:
\begin{equation}
\label{fil-eq}
\nu_t=y_t-\int_0^t E[h_s|\calY_s]ds
\end{equation}
where $(\calY_s, s\in [0,1])$ is the filtration of $y$, called the
observation process.
The following result gives a complete answer to the
 innovation conjecture in the general case to which the above problem
 can be translated:

\begin{theorem}
\label{innov-thm}
 Assume that $U=I_W+u$ is an adapted perturbation of
identity such that $u\in\DD_{2,0}(H)$ and that $E[\rho(-\delta
u)]=1$. Define $L$ as the Radon-Nikodym density 
$$
L=\frac{dU\mu}{d\mu}
$$
and define $v\in L^0_a(\mu,H)$ as $L=\rho(-\delta v)$.
Let $\calU=(\calU_t,t\in [0,1])$ be its filtration
eventually completed with $\mu$-null sets. Let $Z$ be the
innovation process associated to $U$ as defined above, denote by
$\calZ=(\calZ_t,t\in[0,1])$ its filtration. Then $\calU=\calZ$ if
and only if there exists some $\hat{u}\in L^0_a(\mu,H)$ such that
$\hat{U}=I_W+\hat{u}$ is almost surely invertible with inverse
$V=I_W+v$ and  $U=\hat{U}\circ Z$  almost surely.
\end{theorem}
\nproof Sufficiency: We have $\calZ\subset \calU$ by the  construction of $Z$, on
the other hand the relation $U=\hat{U}\circ Z$ implies that
$\calU\subset\calZ$, hence the sufficiency is proved.\\
Necessity:  Suppose now  that $\calZ=\calU$,
let $L$ be the Radon-Nikodym derivative
$$
L=\frac{dU\mu}{d\mu}\,.
$$
Since $L>0$ almost surely, there exists some $v:W\to H$ such that
$\dot{v}$ is adapted and that $L$ can be represented as
$L=\rho(-\delta v)$. Hence the random variable $L$ is represented
by $U$, this implies that $V\circ U=Z$ almost surely, where
$V=I_W+v$. Since $\calU=\calZ$, we can write $U$ as a function of
$Z$, i.e., $U=\hat{U}(Z)$. Then
\beaa
1&=&\mu\{V\circ U=Z\}=\mu\{V\circ \hat{U}(Z)=Z\}\\
&=&\mu\{V\circ\hat{U}(w)=w\}\,,
 \eeaa
 since $Z\mu=\mu$. Consequently, $\hat{U}$ is a right  inverse of
 $V$. Moreover $\hat{U}\mu=\hat{U}\circ Z\mu=U\mu\sim\mu$
 hence it follows from Proposition \ref{inverse-1} that $V\circ
\hat{U}=\hat{U}\circ V=I_W$ $\mu$-almost surely.
 \nqed
\begin{corollary}
Assume that we are in the situation described by the relation 
(\ref{fil-eq}). Let us denote by $\hat{H}:W\to H$ defined by
$$
\hat{H}(t,y)=\int_0^t E[h_s|\calY_s]ds\,.
$$
Denote by $V$ the mapping defined by $V=I_W-\hat{H}$.
Then the filtration generated by the innovation $\nu$
is equal to the filtration of the observation $y$ if and only if 
$$
E\left[\frac{dV}{d\mu}\log \frac{dV\mu}{d\mu}\right]=\frac{1}{2}E[|\hat{H}|_H^2]
$$
\end{corollary}
\nproof
It follows from Theorem \ref{innov-thm}, that the invertibility of $V$
is a necessary and sufficient condition, then we apply Theorem
\ref{entropy-thm}.
\nqed

\begin{remarkk}
In \cite{A-M}, the authors treat the case where the noise is
independent of the signal, this amounts to say that $u$ is
independent of $w$, here on the contrary we are in a situation
where the things are correlated.
\end{remarkk}
\section{The properties of $U\circ V$}
 \noindent
 As we have seen above, the mapping $V\circ U$ preserves
the Wiener measure $\mu$. On the other hand we have, from the
Girsanov theorem \beaa
E[f\circ U\circ V\,L]&=&E[f\circ U\circ V\,\rho(-\delta v)]\\
&=&E[f\circ U]\\
&=&E[f\, L]\,,
\eeaa
for any $f\in C_b(W)$. In other words $U\circ V$ preserves the measure
$\nu$ which is defined by $d\nu=L\,d\mu$. Let us denote  $U\circ V$
with $M$. This mapping is of the form $M=I_W+m$, where $m=v+u\circ V$
is an adapted, $H$-valued mapping.

\begin{proposition}
Assume that $m$ satisfies the following hypothesis:
$$
E[\rho(-\delta m)]=1\,,
$$
where $\delta m$ denotes the It\^o integral of $(\dot{m}_s,s\in
[0,1])$ in $L^0(\mu)$-sense\footnote{This is an abuse of notation
  since the divergence coincides with the It\^o integral only for the
  adapted elements of $L^p(\mu,H)$ with $p>1$.}. Then the mapping
$M=U\circ V$ satisfies the following probabilistic Monge-Amp\`ere equation:
\begin{equation}
\label{MA-2}
L\circ M\,\,E[\rho(-\delta m)|\calM]=E[L|\calM]\,,
\end{equation}
almost surely, where $\calM$ denotes the sigma-algebra generated by $M$.
\end{proposition}
\nproof From the Girsanov theorem, for any $f\in C_b(W)$, we get
$$
E[f\,L]=E[f\circ M\,L\circ M\,\rho(-\delta m)]\,.
$$
On the other hand $M$ preserves the measure $d\nu=L\,d\mu$, hence
$$
E[f\circ M\,L]=E[f\,L]\,.
$$
Therefore
$$
E[f\circ M\,L\circ M\,\rho(-\delta m)]=E[f\circ M\,L]\,,
$$
for any $f\in C_b(W)$ and this proves the claim. \nqed

Let us denote by $(\calM_t,t\in[0,1])$ the filtration generated by $M$
and let us suppose that $m=v+u\circ V$ is in $L^2(\mu,H)$. This last
hypothesis is amply sufficient to ensure the existence of the dual
predictable projection $\hat{m}$ of $m$ with respect to the filtration
$(\calM_t,t\in[0,1])$. It can be calculated  as in Proposition
$$
\hat{m}(t)=\int_0^tE[\dot{m}_s|\calM_s]ds,\,t\in [0,1]\,.
$$
Besides, the innovation  process $(R_t,t\in [0,1])$ associated to $M$,
defined by
$$
R_t=M_t-\int_0^tE[\dot{m}_s|\calM_s]ds
$$
is an $(\calM_t,t\in[0,1])$-Brownian motion and again from \cite{FKK},
any martingale of this filtration can be represented as a stochastic
integral with respect to this innovation process. Consequently, the
martingale $E[\rho(-\delta m)|\calM_t]$ can be represented as in
Proposition \ref{iyi-prop}:
$$
E[\rho(-\delta
m)|\calM_t]=\exp\left(-\int_0^t(E[\dot{m}_s|\calM_s],dR_s)-\frac{1}{2}\int_0^t
|E[\dot{m}_s|\calM_s]|^2ds\right)\,.
$$
From the It\^o representation theorem, there exists an
$(\calM_t,t\in[0,1])$-adapted process $(\dot{\gamma}_t,t\in [0,1])$
such that
$\int_0^1|\dot{\gamma}_t|^2dt<\infty$ almost surely and that
$$
E[L|\calM_t]=\exp\left(-\int_0^t(\dot{\gamma}_s,dR_s)-\frac{1}{2}\int_0^t
|\dot{\gamma}_s|^2ds\right)\,.
$$
Let us calculate  the  terms at the right of  the relation
(\ref{MA-2}):
$$
L\circ M=\exp\left(-\delta v\circ M-\frac{1}{2}|v\circ M|_H^2\right)\,.
$$
Using the identity
$$
\delta v\circ M=\delta(v\circ M)+(v\circ M,m)_H
$$
and taking into account  the exponents of the relation (\ref{MA-2}),
we get
\beaa
&&\delta(v\circ M)+(v\circ M,m)_H+\frac{1}{2}|v\circ M|_H^2\\
&&+\int_0^1(E[\dot{m}_s|\calM_s],dR_s)+\frac{1}{2}\int_0^1|E[\dot{m}_s|\calM_s]|^2ds\\
&=&\int_0^1(\dot{\gamma}_s,dR_s)+\frac{1}{2}|\gamma|_H^2\,,
\eeaa
where the letters without ``dot'' denote the primitives of those with
``dot''. If we restrict all these calculations to the time interval
$[0,t]$, for any $t\in [0,1]$, similar relation holds, consequently we
have proven
\begin{theorem}
If $U\mu=\nu=L\cdot\mu$ and if $L=\rho(-\delta v)$, where $u$ and $v$
are adapted and if $E[\rho(-\delta m)]=1$ and if $m=v+u\circ V\in
L^2(\mu,H)$, then we have the following relation between $v,m$
and $\gamma$:
\begin{equation}
\dot{v}_t\circ M+E[\dot{m}_t|\calM_t]=\dot{\gamma}_t
\end{equation}
$dt\times d\mu$-almost surely, where the scalar product is that of $\R^d$.
\end{theorem}
\section{Relations with the Monge's transport map}
Assume that the density $L$ is in the class $L\log L(\mu)$. It follows
from \cite{F-U2} that there exists an $H-1$-convex element $\varphi$  of
$\DD_{2,1}$ such that the perturbation of identity $T$ defined as
$$
T(w)=w+\nabla\varphi(w)
$$
maps the Wiener measure $\mu$ to $\nu=L\cdot \mu$ and also there is
another map $S=I_W+\nabla\psi$, $\psi\in\DD_{2,1}$ also $H-1$-convex
such that
$$
\mu\left(\{w:\,S\circ T(w)=w\}\right)=1
$$
and
$$
\nu\left(\{w:\,T\circ S(w)=w\}\right)=1\,.
$$
In particular, whenever $\mu$ and $\nu$ are equivalent, then $T$ and
$S$ are inverse to each other $\mu$-almost surely. Let us remark that
neither $T$ nor $S$ are adapted to the filtration $(\calF_t)$. We
shall assume in the sequel that $L$ is $\mu$-almost surely  strictly
positive and represented as before as an exponential density
$L=\rho(-\delta v)$. Let us denote by $(\calT_t, \,t\in [0,1])$ the filtration
generated by $(T_t,\,t\in [0,1])$, where $T_t$ is defined as
$T_t(w)=w(t)+\nabla\varphi(t)$ with
$\nabla\varphi(t)=\int_0^tD_s\varphi ds$. We have
\begin{theorem}
\label{Monge-1}
Assume further that $L\in L^{1+\eps}(\mu)$ for some $\eps>0$, then $T$
is a $\mu$-semimartingale with respect to $(\calT_t)$ and it has the
following decomposition:
\begin{equation}
\label{decom-1}
T_t=B_t+\int_0^t\left(\frac{E[D_s L|\calF_s]}{E[L|\calF_s]}\right)\circ T ds\,,
\end{equation}
where $B=(B_t)$ is a $(\calT_t)$-Brownian motion. Moreover
(\ref{decom-1}) can be also expressed as
\begin{equation}
\label{decom-2}
T_t=B_t-\int_0^t \dot{v}_s\circ T\, ds\,,
\end{equation}
where $\dot{v}$ is defined as $L=\rho(-\delta v)$.
\end{theorem}
\nproof Since $(W_t,t\in  [0,1])$ is the canonical Brownian
motion, the equality $\calT_t=T^{-1}(\calF_t)$ is immediate.
Consequently, for any positive, measurable function $f$, we have
the following identity:
$$
E[f\circ T |\calT_t]=E_\nu[f|\calF_t]\circ T\,.
$$
This relation implies that $(T_t,t\in [0,1])$ is a
$(\mu,(\calT_t))$- quasimartingale if and only if $(W_t,t\in
[0,1])$ is a $(\nu,(\calF_t))$-quasimartingale. This latter
property is immediate since $V=W+v$ is a
$(\nu,(\calF_t))$-Brownian motion and $E_\nu[|v|_H^2]=2E[L\,\log
L]<\infty$. Let us calculate the  drift of $(T_t,\,t\in [0,1])$:
if $\theta$ is a bounded, $\calF_t$-measurable cylindrical
function, we have, using the integration by parts formula \beaa
\frac{1}{h}E[(T_{t+h}-T_t)\theta\circ T]&=&\frac{1}{h}E[(W_{t+h}-W_t)\theta\,L]\\
&=&\frac{1}{h}E\left[\theta\,\int_t^{t+h}D_s L\, ds \right]\\
&\to&E\left[\theta\,D_t L\right]\\
&=&E[\theta E[D_t L|\calF_t]]\\
&=&E\left[\theta E[D_t L|\calF_t]\frac{L}{L_t}\right]\\
&=&E\left[\theta\circ T \frac{E[D_t L|\calF_t]}{L_t}\circ
T\right]\,, \eeaa as $h\to 0$, where $L_s=E[L|\calF_s]$. Moreover,
the local martingale part is a continuous process with $\lip
B^i,B^j\rip_t=\delta_{i,j}t$, hence it is a Brownian motion and
$(T_t)$ has the  decomposition given by the formula
(\ref{decom-1}) which is equivalent to the decomposition given by
(\ref{decom-2}). In fact $L$ can be represented as
$$
L=1+\int_0^1E[D_sL|\calF_s]\,dW_s\,.
$$
On the other hand from the It\^o's formula, we have
$$
L=1-\int_0^1\dot{v}_s\,L_s\, dW_s
$$
hence $L_s\dot{v}_s=-E[ D_sL|\calF_s]$ $ds\times d\mu$-almost
surely. \nqed
\begin{remarkk}
We could have guessed  this theorem by observing simply that the mapping
$B=V\circ T$ preserves the Wiener measure due to the Girsanov theorem.
Therefore the process $(t,w)\to
B(w)(t)$ is a Brownian motion with respect to its own
filtration. However the theorem says that it is also a Brownian motion
with respect to the larger filtration $(\calT_t,t\in [0,1])$.
\end{remarkk}

\begin{theorem}
Assume that $L=\rho(-\delta v)$ satisfies the hypothesis of Theorem
\ref{Monge-1}, let $V=I_W+v$. The map $V$ is not invertible, i.e., the
equation
\begin{equation}
\label{eqn-7}
U_t=W_t-\int_0^t \dot{v}_s\circ U\,ds
\end{equation}
has no strong strong solution if and only if the equation
\begin{equation}
\label{eqn-8}
T_t=B_t-\int_0^t \dot{v}_s\circ T\,ds
\end{equation}
has no strong solution.
\end{theorem}
\nproof Assume that $T$ is a strong solution, then by definition
$T$ should be adapted to the filtration of the Brownian motion
$B=(B_t)$, hence it is of the form $T=\hat{T}\circ B$. Then \beaa
1&=&\mu\{B=\hat{T}\circ B+v\circ \hat{T}\circ B\}\\
&=&\mu\{w=\hat{T}(w)+v\circ \hat{T}(w)\}\\
&=&\mu\{w:\,V\circ \hat{T}(w)=w\}=\mu(D)\,,
\eeaa
hence $\hat{T}$ is a right inverse to $V$. Moreover, for any $f\in
C_b(W)$,
$$
E[f\circ \hat{T}]=E[f\circ \hat{T}\circ B]=E[f\circ T]=E[f\,L]\,.
$$
Therefore $\hat{T}\mu$ is equivalent to $\mu$. Since
$$
1_{\hat{T}(D)}\circ \hat{T}\geq 1_D\,,
$$
we obtain $\mu(\hat{T}(D))=1$ which means that $\hat{T}$ is almost
surely surjective, consequently it is also a left inverse and it
follows from Proposition \ref{inverse-1} that  $\hat{T}$ is a
strong solution to the equation (\ref{eqn-7}), which is a
contradiction. To show the sufficiency suppose that the equation
(\ref{eqn-7}) has a strong solution $U$, then $U$ and $V$ are
inverse to each other almost surely. moreover $B=V\circ T$ is also
invertible hence $U=T\circ B^{-1}$ is $(\calF_t)$-adapted and this
implies that $T$ is $(B^{-1}(\calF_t))$-adapted, consequently the
equation (\ref{eqn-8}) has a strong solution which is a
contradiction. \nqed

\section{Variational techniques for representability and invertibility}
In this section we shall derive a necessary and sufficient condition
for a large class of adapted perturbation of identity. We begin with
some technical results:
\begin{lemma}
\label{var-1-lemma}
Assume that $f\in \DD_{2,1}$ and $\eta \in \DD_{2,0}^a(H)$ such
that $|\eta|_H\in L^\infty(\mu)$. Then we have
$$
f(w+\eta(w))=f(w)+\int_0^1\nabla_\eta f(w+t\eta(w))dt
$$
$\mu$-almost surely.
\end{lemma}
\nproof If $f$ is Fr\'echet differentiable or if it is $H-C^1$,
then the identity is obvious. Assume that $(f_n,n\geq 1)$ is a
sequence of such functions converging to $f$ in $\DD_{2,1}$ and
denote $I_W+\eta$ by $T_\eta$.  Then we have on the one hand \beaa
E[|f_n\circ T_\eta-f_m\circ
T_\eta|]&=&E\left[|f_n-f_m|\,\frac{dT_\eta\mu}{d\mu}\right]\\
&\leq&E[|f_n-f_m|^2]^{1/2}E\left[\left(\frac{dT_\eta\mu}{d\mu}\right)^2\right]^{1/2}\,.
\eeaa
From Lemma \ref{R-N-lemma}, we have
\beaa
E\left[\left(\frac{dT_\eta\mu}{d\mu}\right)^2\right]&=&E\left[\left(\frac{dT_\eta\mu}{d\mu}\right)\circ
T_\eta\right]\\
&=&E\left[\frac{1}{E[\rho(-\delta \eta)|T_\eta]}\right]\\
&\leq&E\left[\frac{1}{\rho(-\delta \eta)}\right]\\
&=&E\left[\exp\left(\delta\eta+\frac{1}{2}|\eta|_H^2\right)\right]<\infty
\eeaa
since $|\eta|_H\in L^\infty(\mu)$. Hence  we get that
$$
\lim_{n,m\to\infty}E[|f_n\circ T_\eta-f_m\circ
T_\eta|]=0\,.
$$
Similarly
\beaa
E\int_0^1|\nabla_\eta f_n-\nabla_\eta
f_m|_H\circ T_{t\eta}dt&=&E\left[|\nabla_\eta f_n-\nabla_\eta
  f_m|_H\int_0^1\frac{dT_{t\eta}\mu}{d\mu}dt\right]\\
&\leq&\|f_n-f_m\|_{2,1}\left(E\int_0^1\left(\frac{dT_{t\eta}\mu}{d\mu}\right)^2\right)^{1/2}\\
&\leq&\|f_n-f_m\|_{2,1}\left(E\int_0^1\exp\left(t\delta\eta+\frac{t^2}{2}|\eta|_H^2\right)dt\right)^{1/2}\\
&&\to 0 \eeaa as $n,m\to \infty$. 
\nqed

\begin{corollary}
Assume that $f\in \DD_{2,1}$ is $\calF_{t_0}$-measurable for some
fixed $t_0<1$. Then the conclusion of Lemma \ref{var-1-lemma} holds
for any $u\in\DD^a_{2,0}(H)$.
\end{corollary}
\nproof
Let $(\tau_n)$ be a sequence of stopping times increasing to infinity
such that $|u^{\tau_n}|$ is essentially bounded where $ u^{\tau_n}$  is defined as 
$$
u^{\tau_n}(t)=\int_0^t 1_{[0,\tau_n]}(s)\dot{u}_sds\,.
$$
From Lemma \ref{var-1-lemma}, it follows trivially that 
$$
f(w+u^{\tau_n}(w))=f(w)+\int_0^1(\nabla f(w+tu^{\tau_n}(w)), u^{\tau_n}(w))_Hdt\,,
$$
moreover, on the set  $\{\tau_n>t_0\}$, we have
$f(w+u^{\tau_n}(w))=f(w+u(w))$ and 
$$
(\nabla f(w+tu^{\tau_n}(w)),
u^{\tau_n}(w))_H=(\nabla f(w+tu(w)), u(w))_H
$$
 almost surely.
\nqed

\noindent

\begin{theorem}
\label{var-f}
 Assume that  $v\in\DD_{2,2}^a(H)$ such that $|v|_H\in L^\infty(\mu)$
 and that 
$$
E[\exp\eps\|\nabla v\|^2_{\rm op}]<\infty
$$
for some $\eps>0$, where $\|\nabla v\|_{\rm op}$ denotes the operator
norm of $\nabla v$. If the following  infimum
$$
 \inf\left(\frac{1}{2}E\left[|\xi+v\circ
   (I_W+\xi)|^2\right]:\,\xi\in\DD_{2,0}^a(H)\right)\,,
 $$
 is attained for some  $u$, then  its value is zero 
and $U=I_W+u$ is inverse of the shift $I_W+v$.
\end{theorem}
\nproof The main point is to show the validity of the variational
formula:
\begin{equation}
\label{id-1}
v(w+u(w)+\eta(w))=v(w+u(w))+\int_0^1\nabla_\eta v(w+u(w)+t\eta(w))dt
\end{equation}
almost surely where $\eta\in\DD_{2,0}^a(H)$ with $|\eta|_H\in
L^\infty(\mu)$  and that these terms are properly integrable in such a
way that the Gateaux derivative at $u$ of $F(u)$ is well-defined. Let
us denote by $v_n$ the regularization of $v$ defined as $P_{1/n}v$,
where $P_{1/n}$ is the Ornstein-Uhlenbeck semigroup. Since $v_n$ is
$H$-differentiable, we get trivially the identity:
\begin{equation}
\label{id-2}
v_n(w+u(w)+\eta(w))=v_n(w+u(w))+\int_0^1\nabla_\eta v_n(w+u(w)+t\eta(w))dt\,.
\end{equation}
By the Jensen inequality we have
\begin{equation}
 \sup_n E\left[\exp\eps\|\nabla v_n\|_{\rm op}\right]<\infty\,\label{exp-ineq}.
\end{equation}
Let us denote by $T_t$ the shift $I_W+u+t\eta$. Then
$$
E\int_0^1|\nabla_\eta v_n\circ T_t|_Hdt\leq\|\eta\|_{L^\infty(\mu)}
E\int_0^1\|\nabla v_n\|_{\rm op}\, l_t\, dt
$$
where $l_t$ is the Radon-Nikodym derivative of $T_t\mu$ with
respect to $\mu$. Using the Young inequality for the dual convex
functions $\exp$ and $x\log x$ we obtain, for any $\kappa>0$,
\begin{equation}
\|\nabla v_n\|_{\rm op}\, l_t\,\leq \exp\kappa\|\nabla v_n\|_{\rm
op}+\frac{1}{\kappa} \,\, l_t\,\log l_t\,\label{young-ineq}.
\end{equation}
It is clear that, from the hypothesis and the Jensen lemma, the sequence
$(\exp\kappa\|\nabla v_n\|,\,n\geq 1)$ is uniformly integrable for
small $\kappa>0$. 
From Lemma \ref{R-N-lemma}
$$
l_t\circ T_t\,E[\rho(-\delta(u+t\eta))|T_t]\leq 1\,,
$$
hence
\beaa
E[l_t\log l_t]&=&E[\log l_t\circ T_t]\\
&\leq& E[-\log E[\rho(-\delta(u+t\eta))|T_t]]\\
&\leq& E[-\log \rho(-\delta(u+t\eta))]\\
&=&\frac{1}{2}E[|u+t\eta|_H^2]\\
&\leq& E[|u|_H^2]+E[|\eta|_H^2]\,. \eeaa 
\noindent 
Hence $(l_t,\,t\in[0,1])$ is uniformly integrable, but we also  need to
prove the uniform integrability of $(l_t\log l_t,t\in[0,1])$. For
this, let $A$ be any measurable subset of $W$, we have, again from  Lemma \ref{R-N-lemma},
\beaa
E[1_Al_t\log l_t]&=&E[1_A\circ T_t\,\log l_t\circ T_t]\\
&=&E[1_A\circ T_t\,(-\log E[\rho(-\delta(u+t\eta))| T_t])]\\
&\leq&E[1_A\circ T_t(\delta(u+t\eta)+\frac{1}{2}|u+t\eta|_H^2)]\\
&\leq&E[1_A\circ T_t\delta(u+t\eta)]+E[1_A\circ
T_t\frac{1}{2}|u+t\eta|_H^2]\,.
\eeaa
The last  two terms are equivalent, hence it suffices to show that the
second terms can be chosen arbitrarily small by choosing $\mu(A)$
small enough. However this is obvious from the integrability of
$|u|_H^2$ and from the uniform integrability of $(l_t,t\in [0,1])$.
From this and
from the inequality (\ref{exp-ineq}), we see that the left hand
side of (\ref{young-ineq}) is uniformly integrable. Consequently
we can pass to the limit in  the relation (\ref{id-2})  in
$L^1(\mu)$ and obtain the relation (\ref{id-1}). We 
can now  calculate the Gateaux derivative of $F$ at $u$ in any
direction $\eta\in \DD_{2,0}^a(H)$ with $|\eta|_H\in
L^\infty(\mu)$ (instead of $\eta\circ U$)  as follows:
\begin{equation}
\label{vf}
F(u+\la\eta)-F(u)=E\int_0^\la
\left(u+t\eta+v\circ(I_W+u+t\eta),(I_H+\nabla
  v)\circ(I_W+u+t\eta)[\eta]\right)_Hdt\,.
\end{equation}
Let us remark that 
\begin{eqnarray}
\lefteqn{E\left[|u|_H\|\nabla v\circ(I_W+u+t\eta)\|_{\rm op}\right]}\nonumber\\
&\leq&E[|u|_H^2]^{1/2}E\left[\|\nabla v\circ(I_W+u+t\eta)\|^2_{\rm
  op}\right]^{1/2}\nonumber\\
&\leq&E[|u|_H^2]^{1/2}E\left[\exp\eps\|\nabla v\|^2_{\rm
  op}+\frac{1}{\eps}l_{t\eta,u}\log l_{t\eta,u}\right]^{1/2}\,,\label{control}
\end{eqnarray}
where 
$$
l_{t\eta,u}=\frac{d(I_W+u+t\eta)\mu}{d\mu}
$$
and from Lemma \ref{R-N-lemma}, we know that 
$$
E[l_{t\eta,u}\log l_{t\eta,u}]\leq \frac{1}{2}E[|u+t\eta|_H^2]\,.
$$
Hence we can commute the expectation with the Lebesgue integral in the
formula (\ref{vf}). Let us denote the expectation of the integrand of
(\ref{vf}) by $F'(u+t\eta)[\eta]$. Since $v\in\DD_{2,2}^a(H)$, using
the formula (\ref{id-1}) for $\nabla v$ instead of $v$ and the
inequality (\ref{control}), we see that
the map $t\to F'(u+t\eta)[\eta]$ is continuous on $[0,1]$. Since $u$
is minimal, we should have $F'(u)[\eta]\geq 0$ for any $\eta$ as
above. Writing the things explicitly:
\beaa
F'(u)[\eta]&=& E\left[\left(u+v\circ U,\left(I_H+\nabla v\circ
      U\right)\eta\right)_H\right]\\
&=&E\left[\left((I_H+\nabla v\circ U)^\star(u+v\circ U),\eta\right)_H\right]\\
&\geq& 0\,.
\eeaa
By the invertibility of $I_H+\nabla v$, we get
$$
u+v\circ U=0
$$
almost surely and this is equivalent to the fact that $U=I_W+u$
and $V=I_W+v$ are inverse to each other. In particular $F(u)=0$.
\nqed

\noindent
 As an application of these kind of variational
calculations in relation with the representability, consider the
problem of calculation of
$$
\inf\left(E\left[\half |\alpha|_H^2+f\circ
(I_W+\alpha)\right]\,:\,\alpha\in\DD_{2,0}^a(H)\right)\,,
$$
where $f:W\to\R$ is a fixed Wiener functional. In fact, as it is
shown in \cite{B-D}, this infimum is equal to $-\log E[\exp-f]$
which is also equal to
\begin{equation}
 \label{meas-variation}
 \inf\left(\int_W
fd\gamma+\int_W\frac{d\gamma}{d\mu}\log\frac{d\gamma}{d\mu}d\mu\right)
\end{equation}
where the infimum is taken w.r.to all the probability measures on
$(W,\calB(W))$ and the latter is uniquely attained at
$$
d\gamma_0=\frac{1}{\int e^{-f}d\mu}e^{-f}d\mu\,.
$$
 In the next theorem we shall give sufficient conditions under
which it is attained:

\begin{theorem}
 \label{variation-thm}
 Assume that $f\in\DD_{2,1}$ is a $1$-convex, bounded
Wiener functional such that
$$
E[\exp\eps|\nabla f|_H]<\infty\,,
$$
for some $\eps>0$. Then the infimum
$$
\inf\left(E\left[\half |\alpha|_H^2+f\circ
(I_W+\alpha)\right]\,:\,\alpha\in\DD_{2,0}^a(H)\right)
$$
is attained at some $u\in\DD_{2,0}^a(H)$ and this adapted vector
field satisfies the following relation:
$$
\dot{u}_t+E[D_tf\circ U|\calF_t]=0
$$
$dt\times d\mu$-almost surely, where $U=I_W+u$. Besides we have
\begin{enumerate}
\item
$$
\frac{dU\mu}{d\mu}=\exp\left(-\int_0^1E_{U\mu}[D_tf|\calF_t]dW_t-
\half\int_0^1|E_{U\mu}[D_tf|\calF_t]|^2dt\right)\,,
$$
where $E_{U\mu}$ denotes the expectation with respect to the
measure $U\mu$, i.e., the image of $\mu$ under $U$.
 \item Let $\dot{v}_t=E_{U\mu}[D_tf|\calF_t]$, denote by $Z$ the
 innovation process associated to $U$, i.e., $Z_t=U_t-\int_0^tE[\dot{u}_s|\calU_s]ds$,
 and define $l$ as
 $$
 l=\exp\left(-\int_0^1E[\dot{u}_t|\calU_t]dZ_t-
 \half\int_0^1|E[\dot{u}_t|\calU_t]|^2dt\right)\,,
 $$
where $\calU_t$ is the sigma algebra
$U^{-1}(\calF_t)=\sigma(W_s+u(s),\,s\leq t)$. Then $E[l]=1$ and we
have
 $$
 l\,\frac{dU\mu}{d\mu}\circ U=l\,\rho(-\delta v)\circ U=1
 $$
 almost surely.
 \end{enumerate}
\end{theorem}
\nproof Let $J(\alpha)$ the expectation above  without $\inf$. For
$\la>0$, let $D_\la=\{\alpha\in \DD_{2,0}^a(H):\,J(\alpha)\leq
\la\}$. Then, for sufficiently large $\la$, $D_\la$ is a
non-empty, convex set. Moreover, if $(\al_n,n\geq 1)\subset D_\la$
converges to some $\al$ in $\DD_{2,0}^a(H)$, then, writing
$A_n=I_W+\al_n$, we have
$$
E\left[\frac{dA_n\mu}{d\mu}\log \frac{dA_n\mu}{d\mu}\right]\leq\half
E[|\al_n|_H^2]\,.
$$
Hence the sequence of Radon-Nikodym densities
$(\frac{dA_n\mu}{d\mu},\,n\geq 1)$ is uniformly integrable. This
property, combined with  Lusin theorem implies that $(f\circ
A_n,n\geq )$ converges to $f\circ A$ in $L^p(\mu)$ for any $p\geq
0$, where $A=I_W+\al$. Therefore $D_\la$ is closed, since it is
convex, it is also weakly closed in $\DD_{2,0}^a(H)$. This implies
that $\al\to J(\al)$ is weakly lower semi continuous (l.s.c.).
Since $D_\la$ is weakly compact, $J$ attains its infimum on
$D_\la$ and the convexity of $J$ implies that this infimum is a
global one. The scalar version of Proposition \ref{var-f} implies
that
 \beaa
0&=& E\left[(u,\al)_H+(\nabla f\circ U,\al)_H\right]\\
&=& E\left[(u,\al)_H+(\pi(\nabla f\circ U),\al)_H\right]\,,
 \eeaa
for any bounded $\al\in \DD_{2,0}^a(H)$, where $\pi$ denotes the
dual predictable projection. Hence we get
$$
\dot{u}_t+E[D_tf\circ U|\calF_t]=0
$$
$dt\times d\mu$-almost surely. Taking the conditional expectation
of this relation with respect to $\calU_t$, we obtain immediately
\begin{equation}
\label{innov-eqn}
  E[\dot{u}_t|\calU_t]+E_{U\mu}[D_tf|\calF_t]\circ U=0
\end{equation}
$dt\times d\mu$-almost surely. It is a simple calculation to see
that the equation (\ref{innov-eqn}) implies
$$
l\,\rho(-\delta v)\circ U=1
$$
almost surely. From the Girsanov theorem, we get
$$
1=E[l\,\rho(-\delta v)\circ U]\leq E[\rho(-\delta v)]\,,
$$
therefore $E[\rho(-\delta v)]=1$. Similarly, for any positive,
measurable $g$ on $W$, we have
$$
E[g\circ U]=E[g\circ U\,l\,\rho(-\delta v)\circ U]\leq
E[g\,\rho(-\delta v)]\,,
$$
therefore
$$
\frac{dU\mu}{d\mu}\leq \rho(-\delta v)\,,
$$
since both are probability densities, they are equal $\mu$-almost
surely. To prove $E[l]=1$ it suffices to write
 $l=1/\rho(-\delta v)\circ U$, then
 \beaa
 E[l]&=&E\left[\frac{1}{\rho(-\delta v)}\circ
U\right]\\
  &=&E\left[\rho(-\delta v)\,.\frac{1}{\rho(-\delta v)}\right]\\
&=&1\,
 \eeaa
 and this completes the proof.
 \nqed

\begin{remarkk}
Suppose that $\|\nabla^2 f\|_{\rm op}\leq c<1$ almost surely,
where $c>0$ is a fixed constant and the norm is the operator norm
on $H$. Then the map $\Phi:\DD_{2,0}^a(H)\to \DD_{2,0}^a(H)$
defined by
$$
\Phi(\xi)=-\pi(\nabla f\circ (I_W+\xi))\,,
$$
where $\pi$ denotes the dual predictable projection, is a strict
contraction, hence there exists a unique $u\in\DD_{2,0}^a(H)$
which satisfies the equation
$$
\dot{u}_t+E[D_tf\circ U|\calF_t]=0
$$
$dt\times d\mu$-almost surely.
\end{remarkk}

\begin{corollary}
Let $u\in\DD_{2,0}^a(H)$ be a minimizer whose existence is assured
by of Theorem \ref{variation-thm}. Define  $U=I_W+u$.  Then
$$
\frac{dU\mu}{d\mu}=\frac{e^-f}{E[e^-f]}=L
$$
{\bf{if and only if}} $U$ is a.s. invertible.
\end{corollary}
\nproof Since
$$
J(u)=E[f\,L]+E[L\log L]=E[f\circ U]+\half E[|u|_H^2]
$$
and since by the hypothesis we have $E[f\,L]=E[f\circ U]$, we
obtain
$$
E[L\log L]=\half E[|u|_H^2]\,.
$$
On the other hand, from Theorem \ref{variation-thm},
 \beaa
  E[L\log L]&=&E[\log L\circ U]\\
  &=&E[-\log l]\\
  &=&\half E\left[\int_0^1|E[\dot{u}_s|\calU_s]|^2ds\right]\,.
  \eeaa
  Consequently, $\dot{u}_s=E[\dot{u}_s|\calU_s]$
  $ds\times d\mu$-almost surely. This implies that $E[\rho(-\delta
  u)]=1$, hence the hypothesis of Theorem \ref{entropy-thm} is
  satisfied and the invertibility of $U$ follows. Conversely,
  suppose that $U$ is invertible, let $M$ be the Radon-Nikodym
  density of $U\mu$ w.r. ro $\mu$. Then we have
  $$
  J(u)=\int_W f\,Md\mu+\int_W M\log M d\mu\,,
  $$
  hence $Md\mu=Ld\mu$ by the uniqueness of the solution of the
  minimization problem (\ref{meas-variation}).
  \nqed

\noindent
{\bf{Acknowledgment:}}This work has been done during my  sabbatical visit to
  the Departement of Mathematics of Bilkent University, Ankara, Turkey.

\vspace{2cm}

{\footnotesize{\bf{
\noindent
A.S. \"Ust\"unel, Telecom-Paristech (formerly ENST),
 Dept. Infres,\\
46, rue Barrault, 75013 Paris, France\\
and\\
Bilkent University, Dept. Math., Ankara, Turkey.\\
email: ustunel@telecom-paristech.fr}
}}

\end{document}